\documentclass[aps,twocolumn,preprintnumbers,shopwacs,shweys,amsmath,amssymb,amsbsy]{revtex4}

\usepackage{pgfplots}
\usepackage{bold-extra}

\def\aslash{\raisebox{0.40ex}{\tiny$\setminus$}}
\def\slash{\raisebox{0.40ex}{\tiny/}}
\def\ish{\mathfrak{I}}
\def\jsh{\mathfrak{J}}

\pgfplotsset{compat=1.18}
\begin{document}

\title{GENERALIZATION OF BINET’S FORMULA FOR FIBONACCI-TYPE NUMERIC SEQUENCES THROUGH THE USE OF ARITHMETIC PSEUDO-OPERATORS}

\author{Victor Vizcarra}

\affiliation{Departamento de Física\\Universidade Estadual de Maringá\\
Av. Colombo, 5790, Maringá - Pr, 87020-900, Brasil.}

\date{\today}

\begin{abstract}
{\indent This paper presents an innovative approach to the study of recurrent sequences by introducing the concept of arithmetic pseudo-operators. Unlike conventional operators, these pseudo-operators are pure complex numbers with specific structural properties, allowing for unprecedented operational reformulations. Represented by the symbols “,+,”, “,\slash,” (slash) e “,\aslash,” (aslash), ``$\bot$", ``$\top$" e ``$\dashv$", these operators correspond to rotations in the unit circle of the complex plane and generalize fundamental operations, as seen in the identities “1\ \slash\ 1\ \aslash\ 1 = 0” and ``1\ $\bot$\ 1\ $\top$\ 1\ $\dashv$\ 1 = 0", which exhibit behavior analogous to conventional subtraction “1 - 1 = 0”. Based on this structure, we reformulate Binet’s equations for the Fibonacci and Tribonacci sequences and outline the path for their generalization to the Tetranacci sequence \cite{koshy}. This new perspective not only enhances the understanding of higher-order recurrences but also suggests potential applications in discrete mathematics and computational algebra, expanding the scope of classical algebraic operations.}.
\end{abstract}

\maketitle

\section{Introduction}

The Fibonacci sequence is defined by the recurrence relation $x_{n+2} = x_{n+1} + x_n$, whose explicit solution can be obtained using Binet’s formula: $x_n = (\phi^n - (1 - \phi)^n)/\sqrt{5}$, where $\phi = \displaystyle\frac{1+\sqrt{5}}{2}$ is the golden ratio. Similarly, the Tribonacci and Tetranacci sequences are defined by the recurrences $x_{n+3} = x_{n+2} + x_{n+1} + x_n$ and $x_{n+4} = x_{n+3} + x_{n+2} + x_{n+1} + x_n$, respectively. The objective of this paper is to generalize these recurrence relations, as well as their closed-form expressions in the style of Binet’s formula. To this end, we explore the properties of arithmetic pseudo-operators: "+", "\,\slash\," and "\,\aslash\,, along with other higher-order pseudo-operators, which will be formally defined and applied throughout the text. These operators play a fundamental role in structuring the closed-form solutions of the studied recurrences, enabling a more systematic approach to obtaining the corresponding algebraic expressions.

\section{Pseudo-Operators}

\noindent{}\textbf{Definition I:} Let $n$ be a positive integer. We define the $n$th root of unity as any complex number $z$ that satisfies the equation
\begin{equation}
    z^n - 1 = 0.
    \label{eq.ciclotomica}
\end{equation}
By the fundamental theorem of algebra, Eq. (\ref{eq.ciclotomica}) has exactly $n$ distinct solutions in the set of complex numbers. These solutions are called the roots of unity.${}$

\vspace{.1in}

\noindent{}\textbf{Property I:} Let $R_n$ be the set of all $n$th roots of unity, given by $R_n = \{\exp(i2k\pi/n)\ :\ k = 0,1, \dots, n-1\}$. The sum of all $n$th roots of unity is equal to zero, that is,
\begin{equation}\nonumber
\sum_{k = 0}^{n-1}\exp\left(\frac{i 2k \pi}{n}\right) = \frac{[\exp(i2\pi/n)]^n - 1}{\exp(i2\pi/n) - 1} = 0
\end{equation}
\hfill$\blacksquare$

\vspace{.1in}

\noindent{}\textbf{Definition II:} (\textsc{Associated Complex}) Let $n$ be a positive integer. We define the negative $n$th root of unity as any complex number $z$ that satisfies the equation
\begin{equation}
z^n + 1 = 0.
\label{eq.ciclotomicacomplexo}
\end{equation}

\noindent{}\textbf{Property II:} Let $C_n$ be the set of all $n$th roots of Eq. (\ref{eq.ciclotomicacomplexo}), given by $C_n = \{\exp(i(2k+1)\pi/n)\ :\ k = 0,1, \dots, n-1\}$. The sum of all negative $n$th roots of unity is equal to zero, that is,
\begin{equation}\nonumber
\sum_{k = 0}^{n-1}\exp\left(\frac{i(2k+1)\pi}{n}\right) = \exp(i\pi/n)\frac{[\exp(i2\pi/n)]^n - 1}{\exp(i2\pi/n) - 1} = 0
\end{equation}
${}$\hfill$\blacksquare$

\vspace{.1in}

For the case of $n = 2$, the solution sets are $R_2 = \{1,-1\}$ and $C_2 = \{+i,-i\}$, where $i = \sqrt{-1}$ denotes the imaginary unit.

\vspace{.1in}

For the case of $n = 3$, the solution sets are $R_3 = \{1,(-1 + i\sqrt{3})/2,(-1 - i\sqrt{3})/2\}$ and $C_3 = \{(1 + i\sqrt{3})/2, -1, (1 - i\sqrt{3})/2\}$.\\[.1cm]

For the case of $n = 4$, the solution sets are  $R_4 = \{1,i,-1,-i\}$ and $C_4 = \{(1 + i)\sqrt{2}/2, (-1 + i)\sqrt{2}/2, (-1 - i)\sqrt{2}/2, (1 - i)\sqrt{2}/2\}$.\\[.1cm]

The \textit{pseudo-operators} are the elements of $R_n$ that behave like arithmetic operators, despite being complex numbers. For this reason, we call them \textit{pseudo-operators}. Similarly, we refer to the elements of $C_n$ as \textit{pseudo-complexes}. These are merely labels, but they are convenient, necessary, and, so to speak, creative and powerful. We will replace each of these numbers with symbols, similarly to the representation  $-1 =  \exp(i\pi)$.\\[.1cm]

\noindent{}\textbf{Definition III:} (\textsc{Pseudo-operators}) Let $\slash1\ = (-1 + i\sqrt{3})/2 = \exp(i2\pi/3)$ and $\aslash1\ = (-1 - i\sqrt{3})/2 = \exp(i4\pi/3)$, two arithmetic pseudo-operators for which the following relation holds: $1\ \slash\ 1\ \aslash\ 1\ = 0$, analogous to the arithmetic identity 1 - 1 = 0. The slash symbol, (\,$\slash$\,), does not represent a division operation but rather a symbolic pseudo-operator, playing a role similar to that of the subtraction symbol (\, -\, ).

Similarly, we define the pseudo-operators $\bot1\ = i = \exp(i\pi/2)$, $\top1\ = -i = \exp(i3\pi/2)$ and $\dashv1\ = -1 = \exp(i\pi)$, which satisfy the relation: $1\ \bot\ 1\ \top\ 1 \dashv\ 1\ = 0$.\\[.1cm]

\noindent{}\textbf{Definition IV:} (\textsc{Pseudo-complex numbers}) Seja $\ish\ = (1 + i\sqrt{3})/2 = \exp(i\pi/3)$, the pseudo-complex number for which the identity $\ish^2 = \slash 1$, holds, similar to the relation $i^2 = -1$ for the imaginary unit. Let also  $\jsh\ = (1 + i)\sqrt{2}/2 = \exp(i\pi/4)$, the pseudo-complex number for which the identities $\jsh^2 = \bot 1$ e $\jsh^4 = -1$.\\[.1cm]

From these definitions, we can organize the sets $R_n$ and $C_n$ for different values of $n$:
\begin{equation}
\left.\begin{array}{lrl}
\textrm{For}\ n = 2, & \textrm{the sets are} & R_2 = \{1,-1\},\\[.1cm]
  & \textrm{and} & C_2 = \{i,-i\}.\\[.2cm]
\textrm{For}\ n = 3, & \textrm{the sets are} & R_3 = \{1,\slash\,1,\aslash\,1\},\\[.1cm]
  & \textrm{and} & C_3 = \{\ish,\slash\,\ish,\aslash\,\ish\}.\\[.2cm]
\textrm{For}\ n = 4, & \textrm{the sets are} & R_4 = \{1,\bot1,\top1,\dashv1\},\\[.1cm]
  & \textrm{and} & C_4 = \{\jsh,\bot\jsh,\top\jsh,\dashv\jsh\}.
\end{array}\right\}
\label{eq.relud}
\end{equation}
Thus, we can generalize this structure for larger values of $n$.

The sets $R_n$ and $C_n$ constitute groups in the strict definition, meaning a set $G$ equipped with a binary operation $G\times G \rightarrow G$, defined by $(\mu_1,\mu_2) \mapsto \mu_1\cdot\mu_2$. This implies that if $\mu_1,\mu_2 \in G$, then, $\mu_1\cdot\mu_2 \in G$. For the case $n = 3$, for example, we consider the set $\mu = \{1, \slash1, \aslash1\}$. We can verify that this set satisfies the following properties:
\begin{itemize}
\item[1.]\textit{Associativity:} For all $\mu_1,\mu_2,\mu_3 \in G$, we have $(\mu_1\cdot\mu_2)\cdot\mu_3 = \mu_1\cdot(\mu_2\cdot\mu_3)$\ .
\item[2.]\textit{Identity element:} There exists an element $e \in G$ such that for every $\mu \in G$, we have $\mu\cdot e = \mu = e\cdot\mu$.
\item[3.]\textit{Inverse element:} For each $\mu \in G$, there exists $\nu \in G$ such that $\mu\cdot\nu = e = \nu\cdot\mu$.
\end{itemize}
These properties can be summarized in the following multiplication table:
\begin{equation}
\fbox{$
\begin{array}{c|ccc}
\cdot & +1 & \slash\, 1 & \aslash\, 1\\\hline
+1 & + 1 & \slash\, 1 & \aslash\, 1\\
\slash\, 1 & \slash\, 1 & \aslash\, 1 & + 1\\
\aslash\, 1 & \aslash\, 1 & + 1 & \slash\, 1\\
\end{array}$}
\label{eq.produto}
\end{equation}

The pseudo-complex number $\ish$ also forms a group, as demonstrated in the following multiplication table:
\begin{equation}
\fbox{$
\begin{array}{c|cccccc}
\cdot & +1 & \slash\,1 & \aslash\,1 & +\ish & \slash\,\ish & \aslash\,\ish\\\hline
+1 & +1 & \slash\,1 & \aslash\,1 & +\ish & \slash\,\ish & \aslash\,\ish\\
\slash\,1 & \slash\,1 & \aslash\,1 & +1 & \slash\,\ish & \aslash\,\ish & +\ish\\
\aslash\,1 & \aslash\,1 & +1 & \slash\,1 & \aslash\,\ish & +\ish & \slash\,\ish\\
+\ish & +\ish & \slash\,\ish & \aslash\,\ish & \slash\,1 & \aslash\,1 & +1\\
\slash\,\ish & \slash\,\ish & \aslash\,\ish & +\ish & \aslash\,1 & +1 & \slash\,1\\
\aslash\,\ish & \aslash\,\ish & +\ish & \slash\,\ish & 1 & \slash\,1 & \aslash\,1\\
\end{array}$}
\label{eq.tabela_complexo_slash}
\end{equation}

\noindent{}This table confirms that $C_3 = \{\ish,\slash\ish,\aslash\ish\}$ forms a group under the multiplication operation defined earlier.

Thus, we obtain the following fundamental properties::
\begin{equation}\nonumber
(\slash\,1)^3 = (\aslash\,1)^3 = +1,\ \ \ \ (\ish)^3 = -1.
\end{equation}

Similarly, $R_4$ e $C_4$ form groups, whose properties can be summarized in the following multiplication tables:
\begin{equation}
\fbox{$
\begin{array}{c|cccc}
\cdot & +1 & \top\, 1 & \bot\, 1 & \dashv\, 1\\\hline
+1 & + 1 & \top\, 1 & \bot\, 1 & \dashv\, 1\\
\top\, 1 & \top\, 1 & \dashv\, 1 & + 1 & \bot\, 1\\
\bot\, 1 & \bot\, 1 & +1 & \dashv\, 1 & \top\, 1\\
\dashv\, 1 & \dashv\, 1 & \bot\, 1 & \top\, 1 & +1\\
\end{array}$}
\label{eq.produto4}
\end{equation}
and,
\begin{equation}
\fbox{$
\begin{array}{c|cccccccc}
\cdot & +1 & \top\, 1 & \bot\, 1 & \dashv\, 1 & +\jsh & \top\, \jsh & \bot\, \jsh & \dashv\, \jsh\\\hline
+1 & +1 & \top\,1 & \bot\,1 & \dashv\, 1 & +\jsh & \top\,\jsh  & \bot\,\jsh & \dashv\, \jsh\\
\top\, 1 & \top\, 1 & \dashv\, 1 & +1 & \bot\, 1 & \top\, \jsh & \dashv\, \jsh & +\jsh & \bot\, \jsh\\
\bot\, 1 & \bot\, 1 & +1 & \dashv\, 1 & \top\, 1 & \bot\, \jsh & +\ish & \dashv\, \ish & \top\,\ish\\
\dashv\,1 & \dashv\, 1 & \bot\, 1 & \top\, 1 & +1 & \dashv\, \jsh & \bot\, \jsh & \top\, \jsh & +\ish\\
+\jsh & +\jsh & \top\, \jsh & \bot\, \jsh & \dashv\,\jsh & \bot\, 1 & +1  & \dashv\, 1 & \top\, 1\\
\top\,\jsh & \top\,\jsh  & \dashv\, \jsh & +\jsh & \bot\, \jsh & +1 & \top\, 1 & \bot\, 1 & \dashv\, 1\\
\bot\,\jsh & \bot\, \jsh & +\jsh & \dashv\, \jsh & \top\,\ish  & \dashv\, 1 & \bot\,1 & \top\, 1 & +1\\
\dashv\,\jsh & \dashv\,\jsh & \bot\, \jsh & \top\, \jsh  & +\jsh & \top\, 1 & \dashv\, 1 & +1 & \bot\,1\\
\end{array}$}
\label{eq.tabela_complexo_dash4}
\end{equation}

It is evident that both pseudo-operators and pseudo-complex numbers constitute cyclic groups.  Explicitly, we have:
\begin{equation}
\left\{\begin{array}{l}
\textrm{Cyclic group generated by}\ \slash1:\\
\mathcal{C}_3 = \langle\, \slash\,1\, \rangle = \{1,\slash\,1,\aslash\,1\}\\[.2cm]
\textrm{Cyclic group generated by}\ \ish:\\
\mathcal{C}_6 = \langle\, \ish\, \rangle = \{1,\slash\,1,\aslash\,1,\ish,\slash\,\ish,\aslash\,\ish\}\\
\end{array}\right.
\end{equation}

\begin{equation}
\left\{\begin{array}{l}
\textrm{Cyclic group generated by}\ \bot1:\\
\mathcal{C}_4 = \langle\, \bot\,1\, \rangle = \{1,\bot1,\dashv1,\top1\} \\[.2cm]
\textrm{Cyclic group generated by}\ \jsh:\\
\mathcal{C}_8 = \langle\, \jsh\, \rangle = \{1,\bot1,\dashv1,\top1,\jsh,\bot\jsh,\dashv\jsh,\top\jsh\}\\
\end{array}\right.
\end{equation}
${}$\\

This approach highlights the algebraic structure of the sets $R_n$ and $C_n$, emphasizing their symmetry and periodicity.\\[.1cm]

\noindent{}\textbf{Additivity:} The additive operation naturally follows consecutive rotations in the complex plane. For example, the operation $2 - 3$ can be interpreted as follows: \textit{Advance 2 units in the direction $\exp(i 0 \pi)$, followed by 3 units in the direction $\exp(i 1 \pi)$}. Similarly, the operation $2\ \slash\ 3$ represents: \textit{Advance 2 units in the direction $\exp(i 0 \pi)$, followed by 3 units in the direction $\exp(i 2\pi/3)$}. The following figure illustrates the operation $2\ \slash\ 3$, where the initial displacement occurs along the positive real axis, followed by a rotation of $2\pi/3$ before the second displacement.
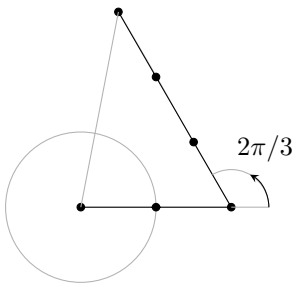
\begin{center}
%\pgfmathsetmacro{\unidade}{1}
%\pgfmathsetmacro{\primervalor}{2}
%\pgfmathsetmacro{\segundovalor}{3}
\begin{tikzpicture}
\draw[color=black!30] (0,0) circle(1);
\draw[fill=black] (0,0) circle(.05cm) -- (1,0) circle(.05cm) -- (2,0) circle(.05cm) -- ++(120:1) circle(.05cm) -- ++(120:1) circle(.05cm) -- ++(120:1) circle(.05cm);
\draw[color=black!30] (2,0) -- (2+.5,0) arc(0:120:.5);
\draw[-stealth] (2,0)  (2+.5,0) arc(0:60:.5);
\draw[color=black!30] (0,0) -- (.5,2.598);
\draw (2,0)  ++(60:.9) node{$2\pi/3$};
\end{tikzpicture}
\end{center}
And just as the operation 2 - 3 has a length, that is, $|2 - 3|$, the operation $2\ \slash\ 3$ also has a well-defined length: $|2\ \slash\ 3| = \sqrt{7}$. Furthermore, the operation $2\ \aslash\ 3$ preserves the same length but with an opposite rotation.\\[.1cm]

\noindent{}\textbf{\textsc{Theorem I:}} Let $R_n$, be a set of arithmetic pseudo-operators associated with rotations in the complex plane. For any $a,b \in \mathbb{R}$, $\mu_1,\mu_2 \in R_n$, the additive operation $\mu_1a + \mu_2b$ has a modulus given by
\begin{equation}
|a\exp(i\eta_1) + b\exp(i\eta_2)| = \sqrt{a^2 + b^2 + 2ab\cos(\eta_1 - \eta_2)}.
\label{eq.comprimento}
\end{equation}
and a resulting angle
\begin{equation}
\theta = \displaystyle\arctan\left(\frac{a\cdot\sin(\eta_1) + b\cdot\sin(\eta_2)}{a\cdot\cos(\eta_1) + b\cdot\cos(\eta_2)}\right)
\label{eq.inclina}
\end{equation}

\noindent{}\textsc{Proof:} Let $\mu \in R_n$ be an arithmetic pseudo-operator represented by a rotation in the complex plane, such that its angles are determined $\eta$, that is, $\mu = \exp(i\eta)$. The elements of $R_n$ satisfy Properties I and II, corresponding to rotational transformations. Consequently, an arithmetic operation such as $\mu_1a + \mu_2b$, analogous to $(-a) + (+b)$ in ordinary arithmetic, can be expressed as the vector sum of complex numbers: $a\exp(i\eta_1) + b\exp(i\eta_2)$, where $\eta_1, \eta_2 \in \eta$. The modulus of this sum corresponds to the length of the resultant vector formed by $a\exp(i\eta_1)$ and $b\exp(i\eta_2)$. Applying the Euclidean norm identity for vector summation in the plane, we obtain equation (\ref{eq.comprimento}). On the other hand, the resulting angle is given by the argument of the obtained complex number. Since the Cartesian projection of the vectors is determined by their real and imaginary components, the tangent of the resulting angle directly follows from their ratio, leading to equation (\ref{eq.inclina}). Therefore, the theorem is proved.\hfill$\blacksquare$

\noindent{}\textbf{\textsc{Example:}} Consider the operation  $2\ \slash\ 3$. According to Property I, the associated set of angles is: $\eta = \{0, \pi/3, 2\pi/3, \pi, 4\pi/3, 5\pi/3\}$. By relation (\ref{eq.relud}), these values correspond to rotations in the complex plane, where each pseudo-operator can be interpreted as a unitary rotation factor. In this case, we define the operation parameters as $a = 2$, $b = 3$, $\eta_1 = 0$, $\eta_2 = 2\pi/3$. Applying equation (\ref{eq.comprimento}) to determine the modulus of the vector sum: $|2\ \slash\ 3| = \sqrt{2^2 + 3^2 + 2\cdot2\cdot3\cdot\cos(0 - 2\pi/3)} = \sqrt{7}$. Now, we determine the resulting angle using equation (\ref{eq.inclina}): $\theta$: $\theta = \arctan\left(\displaystyle\frac{2\cdot\sin(0) + 3\cdot\sin(2\pi/3)}{2\cdot\cos(0) + 3\cdot\cos(2\pi/3)}\right) = \arctan\left(3\sqrt{3}\right)$

The expressions (\ref{eq.comprimento}) and (\ref{eq.inclina}) are general and applicable to any order $n$, with each set $R_n$ having its own pseudo-operators, determined by their respective rotations in the complex plane.

\section{Permutations}

The symbols created for pseudo-operators in the sets $R_3$ and $R_4$ can be used to determine the roots of polynomial equations. We begin with the well-known case of the set $R_2 = \{+,-\}$. A polynomial of degree $n = 2$, expressed in terms of its roots, can be written as: $P(x) = (x - r)(x - s) = x^2 - c_1x - c_0$. Here, we adopt the negative sign in the polynomial form for convenience. The associated elementary symmetric polynomials are given by: $r + s = c_1$ and $r s = -c_0$. We can rewrite $r$ and $s$ as follows: $r = [(r + s) + (r - s)]/2$ and $s = [(s + r) + (s - r)]/2$. If we define $\sigma_1 = r + s$ and $\sigma_2 = s - r = -\sigma_1$, we obtain $r = (c_1 + \sigma_1)/2$ and $s = (c_1 - \sigma_1)/2$. Now, using the relation $rs = -c_0$ from the elementary symmetric polynomials, we find that $\sigma_1 = \sqrt{c_1^2 - 4c_0}$. This leads us to the classical quadratic root formula. However, we will reexpress this result in terms of permutations, thereby establishing the foundation for a more elegant formulation of roots using pseudo-operators.\\

\noindent{}\textbf{Representation in Terms of Permutations:}\\[.2cm]

We can organize this structure in a table that represents the permutations associated with the roots $r$ and $s$:
\begin{center}
\begin{tabular}{c}
(+,+)\\\hline
$r$\ \ $s$\\
$s$\ \ $r$
\end{tabular}\ \ \ \ \ \ \ \ 
\begin{tabular}{c}
(+,\,-)\\\hline
$r$\ \ $s$ \\
$s$\ \ $r$ 
\end{tabular}
\end{center}
In the first column, we identify the sum of the roots $r + s$ and $s + r$, while the second column represents the differences $r - s$ and $s - r$. Both roots are part of these permutations. We can now label each permutation by replacing $r+s$ and $s - r$ with their corresponding values:
\begin{center}
\begin{tabular}{c}
(+,+)\\\hline
$r$\ \ $s$\\
$s$\ \ $r$
\end{tabular}\ \ \ \ \ \ \ \ 
\begin{tabular}{c}
(+,\,-)\\\hline
$r$\ \ $s$ \\
$s$\ \ $r$ 
\end{tabular}\ \ \ \ \ \ $\Rightarrow$\ \ \ \ \ \
\begin{tabular}{c}
(+,+)\\\hline
$c_1$\\
$c_1$
\end{tabular}\ \ \ \ \ \ \ \ 
\begin{tabular}{c}
(+,\,-)\\\hline
$\sigma_1$ \\
$\sigma_2$
\end{tabular}
\end{center}
Since we know that $\sigma_2 = -\sigma_1$, we can rewrite the roots as:
\begin{equation}\nonumber
\left\{\begin{array}{l}
r = \displaystyle\frac12(c_1 + \sigma_1)\\[.5cm]
s = \displaystyle\frac12(c_1 + \sigma_2)
\end{array}\right.\ \ \ \Rightarrow\ \ \ 
\left\{\begin{array}{l}
r = \displaystyle\frac12(c_1 + \sigma_1)\\[.5cm]
s = \displaystyle\frac12(c_1 - \sigma_1)
\end{array}\right.
\end{equation}

Finally, by applying the relation $rs = -c_0$, we obtain $\sigma_1 = \sqrt{c_1^2 + 4 c_0}$. Thus, the permutation structure provides a more organized view of the roots, clearly highlighting the inherent symmetry of the quadratic equation.

For the case $n = 3$, we begin by expressing the elementary symmetric polynomials associated with the cubic polynomial $P(x) = (x - r)(x - s)(x - t) = x^3 - c_2x^2 - c_1x - c_0$. The fundamental relationships between the roots are: $r + s + t = c_2$, $r s + s t + t r = -c_1$ e $r s t = c_0$.

Now, we organize the permutations into a table:

\begin{center}
\begin{tabular}{c}
(+,+,+)\\\hline
$r$\ \ $s$\ \ $t$\\
$t$\ \ $r$\ \ $s$\\
$s$\ \ $t$\ \ $r$
\end{tabular}
\end{center}
In this case, we do not need to list additional permutations such as  $r\ t\ s$, representing $r + t + s$, since the symmetry of the exchanges $r \rightarrow s$, $s \rightarrow t$ and $t \rightarrow r$ ensures that all represent the same quantity. However, when analyzing other permutations, a new structure emerges:
\begin{center}
\begin{tabular}{c}
(+,\,\slash\,,\,\aslash\,)\\\hline
$r$\ \ $s$\ \ $t$\\
$t$\ \ $r$\ \ $s$\\
$s$\ \ $t$\ \ $r$
\end{tabular}\ \ \ \ \ \ 
\begin{tabular}{c}
(+,\,\aslash\,,\,\slash\,)\\\hline
$r$\ \ $s$\ \ $t$\\
$t$\ \ $r$\ \ $s$\\
$s$\ \ $t$\ \ $r$ 
\end{tabular}
\end{center}

\noindent{}Now, we identify the values of the permutations:

\begin{center}
\begin{tabular}{c}
(+,+,+)\\\hline
$c_2$\\
$c_2$\\
$c_2$
\end{tabular}\ \ \ \ \ \ 
\begin{tabular}{c}
(+,\,\slash\,,\,\aslash\,)\\\hline
$\sigma_1$\\
$\sigma_3$\\
$\sigma_5$
\end{tabular}\ \ \ \ \ \ 
\begin{tabular}{c}
(+,\,\aslash\,,\,\slash\,)\\\hline
$\sigma_2$\\
$\sigma_4$\\
$\sigma_6$
\end{tabular}
\end{center}
where we define $\sigma_1 = r\ \slash\ s\ \aslash\ t$, $\sigma_2 = r\ \aslash\ s\ \slash\ t$, and the remaining relations $\sigma_3 = \slash\sigma_1$, $\sigma_4 = \slash\sigma_2$, $\sigma_5 = \aslash\sigma_1$ e $\sigma_6 = \aslash\sigma_2$.

From this, we express the roots in terms of the permutations:
\begin{equation}\nonumber
\left\{\begin{array}{l}
r = \displaystyle\frac13(c_2 + \sigma_1 + \sigma_2)\\[.5cm]
s = \displaystyle\frac13(c_2 + \sigma_5 + \sigma_4)\\[.5cm]
t = \displaystyle\frac13(c_2 + \sigma_3 + \sigma_6)
\end{array}\right.
\end{equation}

\noindent{}or, equivalently,
\begin{equation}\nonumber
\left\{\begin{array}{l}
r = \displaystyle\frac13(c_2 + \sigma_1 + \sigma_2)\\[.5cm]
s = \displaystyle\frac13(c_2\ \aslash\ \sigma_1\ \slash\ \sigma_2)\\[.5cm]
t = \displaystyle\frac13(c_2\ \slash\ \sigma_1\ \aslash\ \sigma_2)
\end{array}\right.
\end{equation}

Finally, using the relations of the symmetric polynomials $r s + s t + t r = -c_1$ and $r s t = c_0$, together with the product properties in the corresponding table, we obtain:
\begin{equation}\nonumber
\left\{\begin{array}{l}
\sigma_1 = \displaystyle\left[\frac12\left(A + \sqrt{A^2 - 4 B^3}\right)\right]^{1/3}\\[.5cm]
\sigma_2 = \displaystyle\left[\frac12\left(A - \sqrt{A^2 - 4 B^3}\right)\right]^{1/3}
\end{array}\right.
\end{equation}
where we define $A = 2 c_2^3 + 9 c_1 c_2 + 27 c_0$ and $B = c_2^2 + 3 c_1$.

Para o caso $n = 4$, the fundamental relationships between the roots are:: $r + s + t + u = c_3$, $r s + r t + r u + s t + s u + t u = -c_2$, $r s t + r s u + r t u + s t u = c_1$ and $r s t u = -c_0$. We define the set of operators as, $R_4 = \{1, \bot, \top, \dashv\}$ where $\bot, \top, \dashv$ are pseudo-operators that obey the fundamental relation: $1\ \bot\ 1\ \top\ 1\ \dashv\ 1 = 0$.

Next, we organize the permutations of the roots into tables:
\begin{center}
\begin{tabular}{c}
(+,+,+,+)\\\hline
$r$\ \ $s$\ \ $t$\ \ $u$\\
$u$\ \ $r$\ \ $s$\ \ $t$\\
$t$\ \ $u$\ \ $r$\ \ $s$\\
$s$\ \ $t$\ \ $u$\ \ $r$
\end{tabular}
\end{center}

\begin{center}
\begin{tabular}{c}
(+,\,$\bot$\,,\,$\top$\,,\,$\dashv$\,)\\\hline
$r$\ \ $s$\ \ $t$\ \ $u$\\
$u$\ \ $r$\ \ $s$\ \ $t$\\
$t$\ \ $u$\ \ $r$\ \ $s$\\
$s$\ \ $t$\ \ $u$\ \ $r$
\end{tabular}\ \ \ \ \ \ 
\begin{tabular}{c}
(+,\,$\bot$\,,\,$\top$\,,\,$\dashv$\,)\\\hline
$r$\ \ $s$\ \ $u$\ \ $t$\\
$t$\ \ $r$\ \ $s$\ \ $u$\\
$u$\ \ $t$\ \ $r$\ \ $s$\\
$s$\ \ $u$\ \ $t$\ \ $r$
\end{tabular}\ \ \ \ \ \ 
\begin{tabular}{c}
(+,\,$\bot$\,,\,$\top$\,,\,$\dashv$\,)\\\hline
$r$\ \ $t$\ \ $s$\ \ $u$\\
$u$\ \ $r$\ \ $t$\ \ $s$\\
$s$\ \ $u$\ \ $r$\ \ $t$\\
$t$\ \ $s$\ \ $u$\ \ $r$
\end{tabular}\ \ \ \ \ \ 
\begin{tabular}{c}
(+,\,$\bot$\,,\,$\top$\,,\,$\dashv$\,)\\\hline
$r$\ \ $t$\ \ $u$\ \ $s$\\
$s$\ \ $r$\ \ $t$\ \ $u$\\
$u$\ \ $s$\ \ $r$\ \ $t$\\
$t$\ \ $u$\ \ $s$\ \ $r$
\end{tabular}\ \ \ \ \ \ 
\begin{tabular}{c}
(+,\,$\bot$\,,\,$\top$\,,\,$\dashv$\,)\\\hline
$r$\ \ $u$\ \ $s$\ \ $t$\\
$t$\ \ $r$\ \ $u$\ \ $s$\\
$s$\ \ $t$\ \ $r$\ \ $u$\\
$u$\ \ $s$\ \ $t$\ \ $r$
\end{tabular}\ \ \ \ \ \ 
\begin{tabular}{c}
(+,\,$\bot$\,,\,$\top$\,,\,$\dashv$\,)\\\hline
$r$\ \ $u$\ \ $t$\ \ $s$\\
$s$\ \ $r$\ \ $u$\ \ $t$\\
$t$\ \ $s$\ \ $r$\ \ $u$\\
$u$\ \ $t$\ \ $s$\ \ $r$
\end{tabular}
\end{center}
Based on these relationships, we can obtain the expressions for the roots $r$, $s$, $t$ $u$. Deriving these algebraic expressions in terms of their discriminants $\sigma$'s is not difficult; it suffices to use the fundamental relationships between the roots along with the properties of the pseudo-operators, as described in the multiplication table (\ref{eq.produto4}). However, the explicit expressions for the $\sigma$'s are extensive and, within the scope of this study, are not essential. The main focus here is the introduction of pseudo-operators, pseudo-complex numbers, and the concept of additivity.

This process can be extended to larger values of $n$, observing that, for each $n$, we have $n$ pseudo-operators, $n!$ permutations involving these pseudo-operators, and $(n-1)!$ values of $\sigma$. When the $\sigma$'s cannot be expressed in terms of radicals, numerical computation becomes a valid approach.

In the next step, we will use these new operators to generalize Binet’s equation.

\section{Generalized Binet Equation}

\subsection{Relationship with Polynomial Equations}

\noindent{}\textbf{Definition IV:} Let $\{x_k\}_{k \geq 0}$ be a sequence defined by the homogeneous linear recurrence relation of order $n$
\begin{equation}
\sum_{j = 0}^{n} a_j x_{k+j} = 0,
\label{eq.seq.geral}
\end{equation}
where $a_j\ \in\ \mathbb{R}$ are constant coefficients. If the ratio between two consecutive terms of this sequence converges to a fixed value $r$, then
\begin{equation}
\lim_{k \rightarrow \infty}\frac{x_{k+1}}{x_k} = r.
\label{eq.lim}
\end{equation}
We call this number the characteristic ratio of the sequence.\\[.1cm]

\noindent{}\textbf{Definition V:} Let $n$ be the order of a homogeneous linear recurrence equation, given by:
\begin{equation}\nonumber
\sum_{j=0}^{n} a_j x_{k+j} = 0,\ \ \ \ \ \textrm{com}\ a_n \neq 0.
\end{equation}

\noindent{}We call the first $n$ values ($x_0, x_1, x_2, \dots,x_{n-1} $), the \textit{seed} or \textit{initial conditions}, which are necessary to uniquely determine the sequence $\{x_k\}_{k\geq0}$. The coefficients $(a_0, a_1, a_2, \dots, a_{n})$, with $a_n \neq 0$, characterize the recurrence equation.\\[.1cm]

\noindent{}\textbf{\textsc{Theorem II:}} Let $\{x_k\}$ be a sequence defined by the recurrence equation of order $n$, given by
\begin{equation}\nonumber
\sum_{j=0}^{n} a_j x_{k+j} = 0,\ \ \ \ \ \textrm{com}\ a_n \neq 0.
\end{equation}

\noindent{} Then, any term $x_{k+n}$ can be expressed as a linear combination of the $n$ previous terms of the sequence, that is,
\begin{equation}\nonumber
x_{k+n} = \sum_{j=0}^{n-1} c_j x_{k+j},
\end{equation}
where the coefficients are given by $c_j = a_j/a_n$.\\[.1cm]

\noindent{}\textsc{Proof:} From the given recurrence equation, we can isolate the term $x_{k+n}$ by rewriting the sum as:
\begin{equation}\nonumber
a_n x_{k+n} = -\sum_{j=0}^{n-1} a_j x_{k+j},
\end{equation}
Since $a_n \neq 0$, we can divide by $a_n$ and obtain
\begin{equation}\nonumber
x_{k+n} = \sum_{j=0}^{n-1} \left(-\frac{a_j}{a_n}\right)\, x_{k+j},
\end{equation}
Defining $c_j = -a_j/a_n$, we obtain the desired expression.\hfill$\blacksquare$\\[.1cm]

\noindent{}\textbf{\textsc{Theorem III:}} Every recurrence equation, as written in (\ref{eq.seq.geral}), is directly associated with a polynomial of degree $n$, of the form 
\begin{equation}
\sum_{j = 0}^{n} a_j r^j = 0.
\label{eq.polinom.geral}
\end{equation}

\noindent{}\textsc{Proof:} We divide Eq. (\ref{eq.seq.geral}) by $x_k$ and rewrite the terms as
\begin{equation}\nonumber
\sum_{j = 0}^{n} a_j \frac{x_{k+j}}{x_k} = 0.
\end{equation}
Using the relation
\begin{equation}\nonumber
\frac{x_{k+j}}{x_k} = \prod_{i = 0}^{j-1} \left(\frac{x_{k+j-i}}{x_{k+j-i-1}}\right),
\end{equation}
and applying the limit property from Eq. (\ref{eq.lim}), we obtain
\begin{equation}\nonumber
\sum_{j = 0}^{n} a_j r^j = 0.
\end{equation}
This shows that $r$ satisfies the characteristic equation associated with the recurrence.\hfill$\blacksquare$\\[.1cm]

\noindent{}\textbf{\textsc{Lemma I:}} The $n$ oots of the polynomial equation are interconnected through elementary symmetric polynomials, and each of these roots is a possible convergence value of the sequence expressed in Eq. (\ref{eq.seq.geral}).\\[.1cm]

\noindent{}\textsc{Proof:} The characteristic polynomial associated with Eq. (\ref{eq.polinom.geral}) can be factored as
\begin{equation}
\prod_{i = 1}^{n} \left(x - r_j\right) = 0.
\label{eq.poly.raiz.geral}
\end{equation}
Expanding this expression and comparing it with the general form of the characteristic polynomial, we obtain the relations between the roots $r_j$ and the coefficients $c_k$ in the form of elementary symmetric polynomials:
\begin{equation}
\left\{\begin{array}{ccc}
r_1 + r_2 + r_3 + \cdots + r_{n} & = & \ \ c_{n - 1}\\[0.2cm]
r_1 r_2 + r_1 r_3 + \cdots + r_{n - 1} r_{n} & = & -c_{n - 2}\\[0.2cm]
\cdots &  & \\[0.2cm]
r_1 r_2 r_3\cdots r_{n} & = & (-1)^{n - 1} c_0
\end{array}\right.\label{polinomiossimetricos}
\end{equation}

\noindent{}Since Definition IV specifies that the ratio between two consecutive terms of the sequence converges to a fixed value $r$, this value must be one of the roots of the characteristic polynomial. Furthermore, any root $r_j$, can play this role, depending on the initial conditions of the sequence. Thus, all roots $r_j$, can be interpreted as possible asymptotic values for the ratio of the sequence terms.\hfill$\blacksquare$

\subsection{Binet’s Equation}

Binet’s equation expresses a specific term of a recurrence sequence in terms of the roots of the associated characteristic polynomial. We can generalize it by considering that the recurrence relation presented in (\ref{eq.seq.geral}) is linear and homogeneous. Thus, any term of the sequence can be written as a linear combination of the solutions of the associated characteristic equation:
\begin{equation}
x_k = \sum_{j = 1}^{n} w_j r_j^k + w_{n+1}.
\label{eq.linear}
\end{equation}
where the coefficients $w_j$ are determined from the initial conditions of the sequence.

In matrix form, this system can be expressed as:
\begin{equation}\left[
\begin{array}{cccccc}
1 & 1 & 1 & \dots & 1 & 1 \\[0.1cm]
r_1 & r_2 & r_3 & \dots & r_{n} & 1 \\[0.1cm]
r_1^2 & r_2^2 & r_3^2 & \dots & r_{n}^2 & 1 \\[0.1cm]
 \vdots &  &  & \vdots &  & \vdots \\[0.1cm]
r_1^n & r_2^n & r_3^n & \dots & r_{n}^n & 1 \\
\end{array}\right]
\left[
\begin{array}{c}
w_1 \\[0.1cm]
w_2 \\[0.1cm]
w_3 \\[0.1cm]
\vdots \\[0.1cm]
w_{n+1} \\
\end{array}\right] = 
\left[
\begin{array}{c}
x_0 \\[0.1cm]
x_1 \\[0.1cm]
x_2 \\[0.1cm]
\vdots \\[0.1cm]
x_n \\
\end{array}\right]\label{sistemageral}
\end{equation}
Solving this system allows us to determine the coefficients $w_j$, ensuring that Binet’s equation provides a closed-form expression for the terms of the sequence.\\

\noindent{}\textbf{Second-Order Binet’s Equation:}\\[.2cm]
Applying the established principles to the particular case where $n = 2$, we can express the recurrence equation and its associated characteristic polynomial as:
\begin{equation}\nonumber
\left\{
\begin{array}{rcl}
\displaystyle x_{k+2} & = & c_1 x_{k+1} + c_0 x_k \\[.2cm]
\displaystyle x^2 & = & c_1 x + c_0 \\
\end{array}\right.
\end{equation}
where we have two initial conditions $(x_0,x_1)$ and two coefficients $(c_0,c_1)$. The roots of the associated quadratic polynomial, which we call generalized golden numbers, are given by
\begin{equation}\nonumber
(r_1,r_2) = \displaystyle\frac{c_1 \pm \sqrt{c_1^2 + 4 c_0}}2,\end{equation}\\

Substituting this expression into the general form of Binet’s equation given in (\ref{eq.linear}), we obtain:
\begin{equation}\nonumber
x_k = w_1 r_1^k + w_2 r_2^k + w_3.
\end{equation}
The coefficients $w_1$, $w_2$ and $w_3$ can be determined by the following linear system:
\begin{equation}\nonumber
\left[
\begin{array}{ccc}
1 & 1 & 1  \\[0.1cm]
r_1 & r_2 & 1 \\[0.1cm]
r_1^2 & r_2^2 & 1
\end{array}\right]
\left[
\begin{array}{c}
w_1 \\[0.1cm]
w_2 \\[0.1cm]
w_3
\end{array}\right] = 
\left[
\begin{array}{c}
x_0 \\[0.1cm]
x_1 \\[0.1cm]
x_2 
\end{array}\right]
\end{equation}
Solving this system results in the generalized Binet’s equation for second-order recurrences:
\begin{equation}
x_k = \frac{(2x_1 - c_1 x_0)}{2}\cdot\left(\frac{r_1^k - r_2^k}{\sqrt{c_1^2 + 4c_0}}\right) + \frac{x_0}{2}\cdot \left(r_1^k + r_2^k\right).
\label{xn2geral}
\end{equation}${}$
This equation provides a closed-form expression for any term $x_k$ of the second-order recurrence sequence, eliminating the need for iterative calculations.\\

However, another approach can be used to derive Binet's equation with the same elegance. Let us observe that $x_k$ can be expressed as a linear combination of the characteristic root powers in the form:
\begin{equation}
x_k = M_1\cdot(r_1^k + r_2^k) + M_2\cdot(r_1^k - r_2^k).
\label{xn2maisgeral}
\end{equation}
Thus, considering the initial values of the sequence, we obtain the conditions
$$x_0 = 2 M_1,\ \ \ x_1 = c_1\cdot M_1 + \sigma_1\cdot M_2.$$
Solving this system for $M_1$ e $M_2$, we find
\begin{equation}\nonumber
M_1 = \frac{x_0}2,\ \ \  M_2 = \frac{2 x_1 - c_1 x_0}{2\sigma_1},
\end{equation}
which leads us directly to equation (\ref{xn2geral}), demonstrating the equivalence between the methods and offering an alternative interpretation of the solution structure.

Note that the general equation contains the terms $\mathcal{F}_k = \left(\displaystyle\frac{r_1^k - r_2^k}{\sqrt{c_1^2 + 4c_0}}\right)$ and $\mathcal{L}_k = r_1^k + r_2^k$, owhich, in the particular case where the coefficients are $(c_0,c_1) = (1,1)$, correspond, respectively, to the Fibonacci and Lucas sequences. The sequence number represented by $x_k$ in Eq. (\ref{xn2geral}), only when $(c_0,c_1) = (1,1)$, coincides exactly with the Fibonacci sequence if the initial conditions are $(x_0,x_1) = (0,1)$, and with the Lucas sequence if we adopt $(x_0,x_1) = (2,1)$. 

Now, keeping  $(x_0,x_1) = (0,1)$ fixed, we obtain the equation $x_k = \mathcal{F}_k$, but with a golden number different from $\phi = (1+\sqrt{5})/2$, since the generalized golden numbers depend on the coefficients $(c_0,c_1)$. This means that by varying the values of the initial conditions and the coefficients, we can generate a wide variety of sequences, including cases with real and complex numbers, revealing unusual and sometimes counterintuitive patterns.\\[0.1cm]

\noindent{}\textbf{Binet’s Equation of Order 3:}\\[.2cm]
For $n = 3$, the linear recurrence takes the form:
\begin{equation}\nonumber
\left\{\begin{array}{rcl}
\displaystyle x_{k+3} & = & c_2 x_{k+2} + c_1 x_{k+1} + c_0 x_k \\[.2cm]
\displaystyle x^3 & = & c_2 x^2 + c_1 x + c_0 \\
\end{array}\right.
\end{equation}
with three initial seeds $(x_0,x_1,x_2)$ and three coefficients $(c_0,c_1,c_2)$.

For the general solution, we adopt an alternative approach to that presented in (\ref{eq.linear}), using instead the structure proposed in (\ref{xn2maisgeral}). In this way, the general solution can be expressed as:
\begin{equation}
\begin{array}{l}
x_k = M_1\cdot(r_1^k + r_2^k + r_3^k) + M_2\cdot(r_1^k\ \slash\ r_2^k\ \aslash\ r_3^k)\\[.3cm]
\hspace{1.5cm}+ M_3\cdot(r_1^k\ \aslash\ r_2^k\ \slash\ r_3^k),
\end{array}
\label{xn3maisgeral}
\end{equation}
where $r_1$, $r_2$ and $r_3$ are the roots of the associated characteristic equation. From the matrix system and by exploring the structure of pseudo-operators, we obtain an elegant formulation of Binet's generalized equation for third-order recurrences:
\begin{widetext}
\begin{equation}
\begin{array}{c}
x_k = \displaystyle\frac{9\,\sigma_1\,x_2 - 3\,(2\,c_2\,\sigma_1 + \sigma_2^2)\,x_1 - ((c_2^2 + 6\,c_1)\, \sigma_1 - c_2\, \sigma_2^2)\,x_0}{3}\left(\frac{r_1^k\ \aslash\ r_2^k\ \slash\ r_3^k}{\sigma_1^3 - \sigma_2^3}\right)\\[.5cm]
\hspace{2cm}-\,\displaystyle\frac{9\,\sigma_2\,x_2 - 3\,(2\,c_2\,\sigma_2 + \sigma_1^2)\,x_1 - ((c_2^2 + 6\,c_1)\, \sigma_2 - c_2\, \sigma_1^2)\,x_0}{3}\left(\frac{r_1^k\ \slash\ r_2^k\ \aslash\ r_3^k}{\sigma_1^3 - \sigma_2^3}\right)\\[.5cm]
\hspace{-2cm}+\,\displaystyle\frac{x_0}{3}\left(r_1^k + r_2^k + r_3^k\right)
\end{array}
\end{equation}
\end{widetext}
where $\sigma_1^3 - \sigma_2^3 = \sqrt{\left(27\,c_0 + 2\,c_2^3 + 9\,c_1\,c_2\right)^2 - 4\,\left(c_2^2 + 3\,c_1\right)^3}$.

This approach not only unifies the representation of solutions but also reveals new algebraic and structural connections between the equation’s coefficients and the characteristic root powers.

In the particular case where the coefficients are $(c_0,c_1,c_2) = (1,1,1)$ and the seeds are $(x_0,x_1,x_2) = (0,1,1)$ we recover the Tribonacci sequence. On the other hand, by choosing $(x_0,x_1,x_2) = (3,1,3)$, we obtain what we would call the tri-Lucas sequence.

Analogously to the second-order case, the solution $x_k$ contains three fundamental terms: $\mathcal{A}_k = \displaystyle\left(\frac{r_1^k\ \aslash\ r_2^k\ \slash\ r_3^k}{\sigma_1^3 - \sigma_2^3}\right)$, $\mathcal{B}_k = \displaystyle\left(\frac{r_1^k\ \slash\ r_2^k\ \aslash\ r_3^k}{\sigma_1^3 - \sigma_2^3}\right)$ and $\mathcal{C}_k = \left(r_1^k + r_2^k + r_3^k\right)$, os which individually exhibit intriguing properties and deserve further analysis, which will be addressed in a future study.\\[0.1cm]

\noindent{}\textbf{Fourth-Order Binet Equation:}\\[.2cm]
Binet's equations generalize the solution of homogeneous linear recurrences with constant coefficients, allowing the general term of the sequence to be expressed in terms of the roots of the characteristic equation. For the case $n = 4$, the linear recurrence takes the form:
\begin{equation}\nonumber
\left\{\begin{array}{rcl}
\displaystyle x_{k+4} & = & c_3 x_{k+3} + c_2 x_{k+2} + c_1 x_{k+1} + c_0 x_k \\[.2cm]
\displaystyle x^4 & = & c_3 x^3 + c_2 x^2 + c_1 x + c_0 \\
\end{array}\right.
\end{equation}
with four initial seeds $(x_0,x_1,x_2,x_3)$ and four coefficients $(c_0,c_1,c_2,c_3)$. The general solution can be expressed as:
\begin{widetext}
\begin{equation}
\begin{array}{l}
x_k = M_1\cdot(r_1^k + r_2^k + r_3^k + r_4^k) + M_2\cdot(r_1^k\ \bot\ r_2^k\ \top\ r_3^k\ \dashv\ r_4^k)\\[.3cm]
\hspace{2.cm} +\ M_3\cdot(r_1^k\ \dashv\ r_2^k\ \bot\ r_3^k\ \top\ r_4^k) + M_4\cdot(r_1^k\ \top\ r_2^k\ \dashv\ r_3^k\ \bot\ r_4^k)\\[.3cm]
\hspace{3.cm} +\ M_5\cdot(r_1^k\ \bot\ r_2^k\ \dashv\ r_3^k\ \top\ r_4^k) + M_6\cdot(r_1^k\ \top\ r_2^k\ \bot\ r_3^k\ \dashv\ r_4^k)\\[.3cm]
\hspace{4.cm} +\ M_7\cdot(r_1^k\ \dashv\ r_2^k\ \top\ r_3^k\ \bot\ r_4^k),
\end{array}
\label{xn4maisgeral}
\end{equation}
\end{widetext}
where $r_1$, $r_2$, $r_3$, $r_4$ are the roots of the associated characteristic equation.

Binet's equation for second-order recurrences, such as the Fibonacci sequence, is relatively simple and well known. However, as the order $n < 4$, increases, obtaining a closed-form solution follows the same principle, but the algebraic manipulation becomes significantly more laborious.

The greatest challenge in generalizing to higher orders lies in obtaining the expressions for the coefficients $M_j$, which depend on the initial conditions and the roots of the characteristic equation. Although they can be written explicitly, constructing these expressions involves extensive calculations and can quickly become impractical for larger values of $n$.

In this paper, we do not present the explicit algebraic expressions for $M_1, \dots, M_7$, not because they are impossible to determine, but due to their considerable complexity and length. The search for these expressions requires arduous work, with voluminous algebraic manipulations involving matrix determinants and linear systems that grow rapidly in size as $n$ increases.

Thus, while Binet's equations remain a powerful tool in solving homogeneous linear recurrences, their practical application to higher orders requires advanced algebraic techniques and often the aid of computational methods to obtain explicit solutions.

\section{Future Perspectives}

The approach of arithmetic pseudo-operators introduced in this work offers a new framework for analyzing recurrent sequences and structuring algebraic relations through innovative operators. However, even broader possibilities emerge when extending this formulation to higher-dimensional spaces.

A promising direction is the generalization of pseudo-operators from $R_3$ to the unit quaternionic sphere, replacing $i = \sqrt{-1}$ with the operator $\hat{u} = (i + j + k)/\sqrt{3}$, where $i^2 = j^2 = k^2 = ijk = -1$. This would allow a new class of operational transformations analogous to rotations in the complex plane, but now embedded in a higher-dimensional space. This expansion suggests connections with the Lie groups $SU(2)$ and $SU(3)$ \cite{cvitanovic}, which are fundamental in the mathematical description of particle physics \cite{griffiths}.

Furthermore, this approach may have implications in Quantum Chromodynamics (QCD), the theory that describes the strong interaction between quarks and gluons \cite{peskin}. Pseudo-complex operators could be analyzed from the perspective of algebraic representations used in QCD, particularly in the decomposition of color states into representations of the $SU(3)$. group. Since quaternions naturally appear in particle physics and quantum symmetries \cite{Adler}, exploring their extensions within the algebra of pseudo-operators may provide new mathematical perspectives for modeling fundamental interactions.

These research directions, in addition to deepening the theory of pseudo-operators, pave the way for applications in computational algebra, theoretical physics, and even in the modeling of complex dynamical systems.

\section{Conclusion}

In this paper, we presented an innovative approach to the generalization of numerical sequences, introducing the pseudo-operators $\slash$ and $\aslash$ the pseudo-complex field $\ish$, and the possibility of defining new higher-degree pseudo-operators, such as $\bot$, $\top$ and $\dashv$, along with their pseudo-complex field $\jsh$. However, far from being mere operators in the conventional sense, these objects should be understood as numbers endowed with specific structural properties, granting them a distinct operational behavior. This new perspective allows us to perceive numbers not merely as abstract quantities but as dynamic mathematical entities carrying deep algebraic and structural relationships.

The generalized Binet equations formulated in this context demonstrate the elegance and power of this approach, enabling the expression of recurrence equation solutions that encompass classical sequences such as Fibonacci and Tribonacci, as well as more general extensions. The introduction of these new numbers significantly expands the available tools for manipulating roots of unity and pseudo-complex numbers, unveiling new connections within the theory of numerical sequences.

However, this formulation presents theoretical and computational challenges. The geometric interpretation of pseudo-operators, their integration with known algebraic structures, and the search for more efficient representations are issues requiring further investigation. The generalization to higher orders demands a deeper understanding of their algebraic properties and the interactions between their pseudo-complex components.

Another promising avenue is the study of the functions $\mathcal{A}_k$, $\mathcal{B}_k$ and $\mathcal{C}_k$ as potential fibonometric functions, analogous to trigonometric functions, capable of establishing new links between numerical sequences and special functions. Investigating these functions could lead to a more comprehensive formulation of the properties of these sequences and their applications in various areas of mathematics and theoretical physics.

Thus, this paper not only proposed a new formalism for analyzing numerical sequences but also suggested a fresh perspective on numbers as structure-bearing entities. We hope this approach inspires future research and contributes to the advancement of mathematical understanding, opening doors to new applications in pure mathematics, theoretical physics, engineering, and computer science.

\end{document}